\documentclass[11pt]{amsart}
\usepackage{amssymb,latexsym,amsmath,graphicx}

\theoremstyle{plain}
\newtheorem{theorem}{Theorem}
\numberwithin{theorem}{section}

\newtheorem{proposition}[theorem]{Proposition}

\theoremstyle{definition}
\newtheorem{definition}[theorem]{Definition}

\begin{document}

\title{A New Way to Compute the Thurston-Bennequin Number}

\author{Hao Wu}

\address{Department of mathematics and Statistics\\
         Lederle Graduate Research Tower\\
         710 North Pleasant Street\\
         University of Massachusetts\\
         Amherst, MA 01003-9305\\
         USA}

\email{wu@math.umass.edu}

\begin{abstract}
We establish a relation between several easy-to-calculate numerical invariants of generic knots in $\mathbb{H}^2 \times S^1$, and use it to give a new method of computing the Thurston-Bennequin number of Legendrian knots in the tight contact $\mathbb{R}^3$.
\end{abstract}

\maketitle

\section{introduction}

A contact structure $\xi$ on $\mathbb{R}^3$ is an oriented positively nowhere integrable tangent plane distribution, i.e., $\xi$ is oriented as a $2$-plane distribution, and there is a $1$-form $\alpha$ on $\mathbb{R}^3$ such that $\xi=\ker\alpha$, $d\alpha|_{\xi}>0$, and $\alpha\wedge d\alpha>0$. Such a $1$-form is called a contact form for $\xi$. A knot in $\mathbb{R}^3$ is said to be Legendrian in $(\mathbb{R}^3,\xi)$ if it is tangent to $\xi$ everywhere. $\xi$ is said to be overtwisted if there is an embedded disk $D$ in $\mathbb{R}^3$ such that $\partial D$ is Legendrian, but $D$ is transversal to $\xi$ along $\partial D$. A contact structure that is not overtwisted is called tight. In \cite{E2,E3}, Y. Eliashberg proved that any tight contact structure on $\mathbb{R}^3$ is isotopic to the standard contact strcture $\xi_0$ given by the contact form $\alpha_0=dz-ydx$, where $(x,y,z)$ are the standard Cartesian coordinates of $\mathbb{R}^3$.

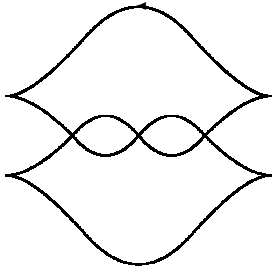
\begin{figure}[h]

\setlength{\unitlength}{1pt}

\begin{picture}(420,100)(-210,-50)

\linethickness{.5pt}

\put(-1,49){\vector(-1,0){0}}

\qbezier(-50,15)(-40,15)(-20,37.5)

\qbezier(-20,37.5)(0,60)(20,37.5)

\qbezier(20,37.5)(40,15)(50,15)

\qbezier(-50,15)(-40,15)(-25,0)

\qbezier(-25,0)(-12.5,-15)(0,0)

\qbezier(0,0)(12.5,15)(25,0)

\qbezier(25,0)(40,-15)(50,-15)

\qbezier(-50,-15)(-40,-15)(-20,-37.5)

\qbezier(-20,-37.5)(0,-60)(20,-37.5)

\qbezier(20,-37.5)(40,-15)(50,-15)

\qbezier(-50,-15)(-40,-15)(-25,0)

\qbezier(-25,0)(-12.5,15)(0,0)

\qbezier(0,0)(12.5,-15)(25,0)

\qbezier(25,0)(40,15)(50,15)

\end{picture}
\caption{Front projection of a Legendrian right hand trefoil knot with $tb=1$, $\mu=0$. Note that at a crossing in a front projection, the branch with smaller slope is always on top.}\label{trefoileg}
\end{figure}

There are two "classical" invariants of Legendrian knots in a contact $\mathbb{R}^3$, the Thurston-Bennequin number $tb$ and the Maslov index $\mu$. We call the framing of a Legendrian knot given by the contact planes its contact framing. The Thurston-Bennequin number is defined to be the index of the contact framing of a Legendrian knot relative to its Seifert framing. And the Maslov index is the obstruction of extending the unit tangent vector of a Legendrian knot to a non-vanishing section of the contact plane bundle over its Seifert surface. See, e.g., \cite{Ben,E,FT} for more about these invariants. In $(\mathbb{R}^3,\xi_0)$, there are two methods to calculate these invariants. First, one can project a Legendrian knot $L$ onto the $xy$-plane, and get an immersed knot diagram, which we call the Legendrian projection. Then the Thurston-Bennequin number is the writhe of the diagram, and the Maslov index is the degree of the Gauss map of the diagram. Second, one can project $L$ onto the $xz$-plane, and get a knot diagram with cusps but no vertical tangents, which we call the front projection. Then the Thurston-Bennequin number is the writhe of the diagram minus half of the number of cusps, and the Maslov index is half of the number of cusps tranversing downward minus half of the number of cusps tranversing upward. See Figure \ref{trefoileg} for an example.

In this paper, we consider the contact structure $\xi_1$ given by $\alpha_1=dz+r^2d\varphi$, where $(z,r,\varphi)$ are the standard cylindrical coordinates of $\mathbb{R}^3$. It is isotopic to the standard contact structure (see, e.g., \cite{Et}). We first define some easy to calculate knot invaraints in $\mathbb{H}^2\times S^1$, and then use these to compute $tb$ and $\mu$ for Legendrian knots in $(\mathbb{R}^3,\xi_1)$ disjoint from the $z$-axis. (This is the generic situation since any Legendrian knot can be make disjoint from the $z$-axis by small Legendrian isotopy.)

\section{Definitions of Invariants}\label{inv-def}

The $2$-dimensional hyperbolic space $\mathbb{H}^2$ is defined to
be the upper half plane $\{(x,y)\in\mathbb{R}^2|~y>0\}$ equipped
with the metric $ds^2=\frac{dx^2+dy^2}{y^2}$. $ST\mathbb{H}^2$ is
the unit tangent bundle of $\mathbb{H}^2$, i.e., the $S^1$-bundle
over $\mathbb{H}^2$ formed by tangent vectors of length $1$. There
is a natural orientation preserving diffeomorphism
$\Phi:\mathbb{H}^2\times S^1 \rightarrow ST\mathbb{H}^2$ given by
\begin{equation}
\Phi(x,y,\theta) ~=~ y(\cos{\theta}\frac{\partial}{\partial x} +
\sin{\theta}\frac{\partial}{\partial y})|_{(x,y)}.
\end{equation}
Let $\mathbb{R}_0^3$ be $\mathbb{R}^3$ with the $z$-axis removed. Then there is a natural orientation preserving
diffeomorphism $\Psi:\mathbb{H}^2\times S^1 \rightarrow
\mathbb{R}_0^3$ given by
\begin{equation}
\left\{%
\begin{array}{ll}
    z ~=~ x,  \\
    r ~=~ y^{\frac{1}{2}}, \\
    \varphi ~=~ \theta. \\
\end{array}%
\right.
\end{equation}
In the rest of this paper, we will always identify these three
manifolds by the diffeomorphisms $\Phi$ and $\Psi$.

An embedding of $S^1$ into $\mathbb{H}^2 \times S^1$ is a called
knot. The standard orientation of $S^1$ induces an orientation on
every knot. Unless otherwise specified, all knots in this paper
are oriented this way. Since $H_1(\mathbb{H}^2 \times S^1)\cong
H_1(S^1)\cong \mathbb{Z}$, and the homology class of a $S^1$-fiber
is a generator of $H_1(\mathbb{H}^2 \times S^1)$, we have that,
for any knot $K$ in $H_1(\mathbb{H}^2 \times S^1)$, there is a
unique integer $h(K)$ satisfying $[K]=h(K)[\{\text{pt}\}\times
S^1]\in H_1(\mathbb{H}^2 \times S^1)$.

\begin{definition}
$h(K)$ is called the homology of the knot $K$.
\end{definition}

A knot in $\mathbb{H}^2 \times S^1$ is said to be generic if it is
nowhere tangent to the $S^1$-fibers. Two generic knots are said to
be generically isotopic if they are isotopic through generic
knots. Let $Pr:\mathbb{H}^2 \times S^1 \rightarrow \mathbb{H}^2$
be the projection onto the first factor. For any generic knot $K$,
$Pr(K)$ is an immersed curve in $\mathbb{H}^2$.

For an immersed curve $L$ in $\mathbb{H}^2$, the canonical lifting
of $L$ to $\mathbb{H}^2\times S^1$ is defined to be
$\widetilde{L}=\frac{dL}{ds}\in ST\mathbb{H}^2 =
\mathbb{H}^2\times S^1$, where $ds$ is the arc length element of
$L$.

\begin{definition}
For a generic knot $K$ in $\mathbb{H}^2\times S^1$, define the
rotation number of $K$ to be $r(K)=h(\widetilde{Pr(K)})$, where
$\widetilde{Pr(K)}$ is the canonical lifting of $Pr(K)$.
\end{definition}

Let $V$ be the tangent vector field on $\mathbb{H}^2\times S^1$
formed by the unit tangent vectors of the $S^1$-fibers. For any
generic knot $K$, $V$ gives $K$ a framing. We call this framing
the canonical framing for $K$. Also, as a knot in
$\mathbb{R}_0^3\subset \mathbb{R}^3$, $K$ has a Seifert surface in
$\mathbb{R}^3$. This surface gives $K$ a framing, called the
Seifert framing for $K$. Note that the Seifert framing is well
defined, i.e., it is independent of the choice of the Seifert
surface.

\begin{definition}
The self linking number $\beta(K)$ of a generic knot $K$ is
defined to be the index of the canonical framing with respect to
the Seifert framing.
\end{definition}

Equivalently, $\beta(K)$ can be defined to be the intersection
number of $K'$ and $\Sigma\cap\mathbb{R}_0^3$, where $K'$ is the
knot given by pushing $K$ slightly in the direction of $V$, and
$\Sigma$ is a Seifert surface of $K$ in $\mathbb{R}^3$.

\begin{figure}[h]
  \includegraphics[width=4in]{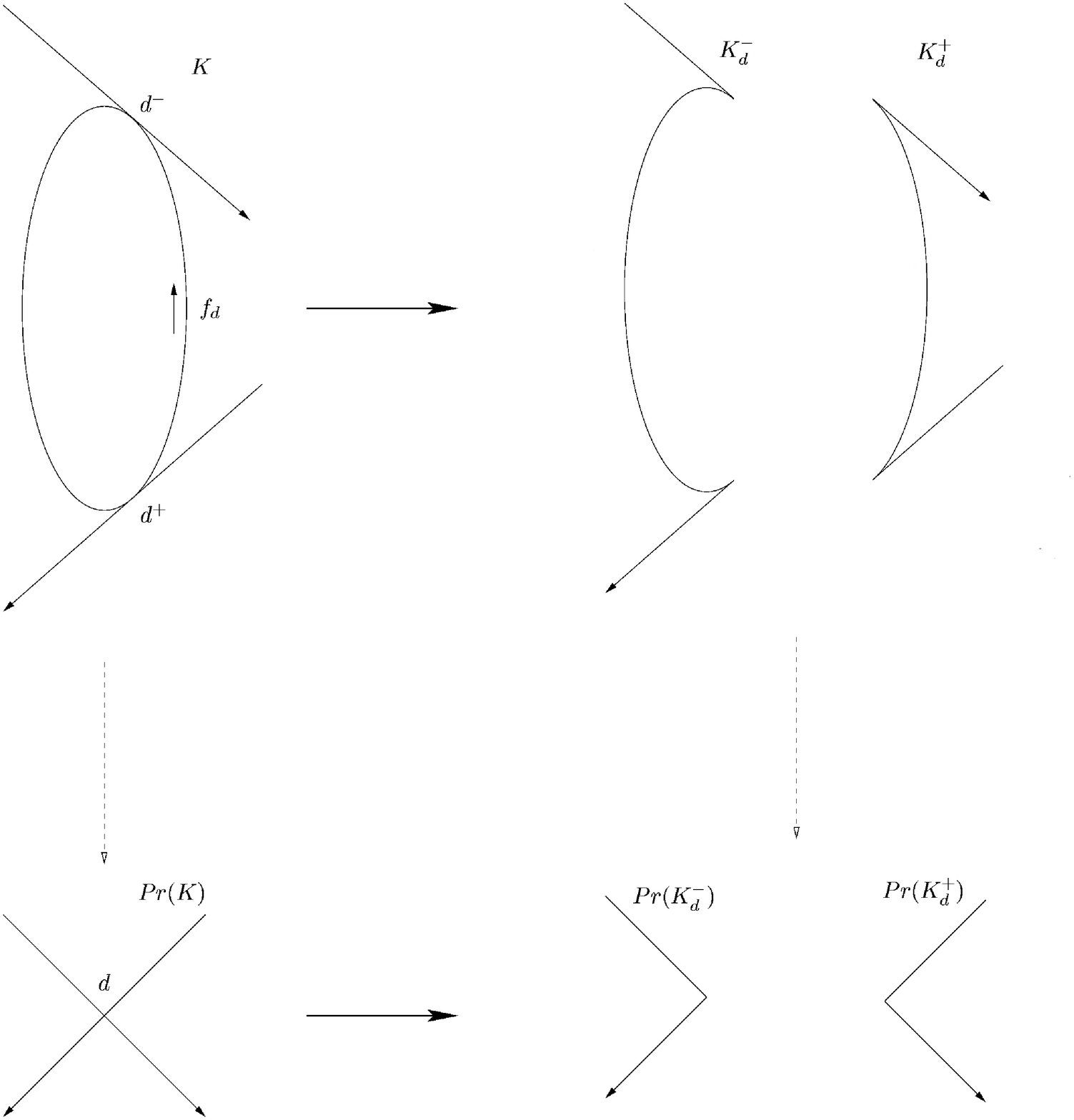}
  \caption{}\label{local}
\end{figure}

Let $K$ be a generic knot, and $d$ a transversal double point of
$Pr(K)$. Let $d^+,d^-\in K$ be the two pre-images of $d$. The
$\pm$ signs here are assigned in such a way that
$\{Pr_*(\frac{dK}{ds}|_{d^+}),Pr_*(\frac{dK}{ds}|_{d^-})\}$ is an
oriented basis for $T_d\mathbb{H}^2$. Note that $d^+$ and $d^-$
are on the same $S^1$-fiber. Denote this fiber by $f_d$. Let
$f_d^+$ be the part of $f_d$ starting from $d^+$ and ending in
$d^-$, and $f_d^-$ the part of $f_d$ starting from $d^-$ and
ending in $d^+$. Cutting $K$ open at $d^+$ and $d^-$, we get two
embedded open oriented curves $C_d^+$ and $C_d^-$ in
$\mathbb{H}^2\times S^1$, where $C_d^+$ is the one starting from
$d^-$ ending in $d^+$, and $C_d^-$ is the one starting from $d^+$
ending in $d^-$. Let $K_d^+$ be the knot formed by $C_d^+\cup
f_d^+$, and $K_d^-$ the knot formed by $C_d^-\cup f_d^-$. Define
the weight $w_d$ of $d$ to be $w_d=h(K_d^+)-h(K_d^-)$.

We can always isotope a generic knot slightly (hence, generically)
so that the only singularities of the projection on $\mathbb{H}^2$
of the perturbed generic knot are transversal double points.

\begin{definition}
Let $K$ be a generic knot, and $K'$ a generic knot generically
isotopic to $K$ such that the only singularities of $Pr(K')$ are
transversal double points. The winding number $w(K)$ of $K$ is
defined to be sum of the weights of all the transversal double
points of $Pr(K')$.
\end{definition}

We need demonstrate that the winding number is well defined and is
invariant under generic isotopy. To do this, we only need show
that, if two generic knots are generically isotopic, and the only
singularities of their projections on $\mathbb{H}^2$ are
transversal double points, then the sums of the weights of all the
transversal double points of the two knots are equal. The proof of
this fact is straightforward. Note that, if two such knots are
generically isotopic, we can change one of them to the other by a
series of local moves corresponding to second type and third type
Reidemeister moves of its projection. There is a natural
correspondence between the transversal double points of two
projections differed by a third type Reidemeister move. One easily
checks that the weights of corresponding double points are the
same. It is also easy to check that, in a second type Reidemeister
move, the sum of the weights of the two transversal double points
to be created or to be cancelled is always $0$. Thus, the winding
number is indeed well defined and invariant under generic isotopy.

\section{Relating $h$, $r$, $\beta$ and $w$}

\begin{theorem}\label{main}
Let $K$ be a generic knot in $\mathbb{H}^2 \times S^1$. Then
\[
\beta(K)=w(K)+r(K)h(K).
\]
\end{theorem}

To prove Theorem \ref{main}, we need use the universal first order
invariant $I$ defined by Chernov in \cite{Cher}. We now briefly
introduce the definition of this invariant in our settings. For
more details, see \cite{Cher}.

Choose a base point $x_0$ in $\mathbb{H}^2 \times S^1$, and let
$\mathcal{B}$ be the group $\pi_1(\mathbb{H}^2 \times S^1,x_0)
\oplus \pi_1(\mathbb{H}^2 \times S^1,x_0)$ modulo the actions of
$\pi_1(\mathbb{H}^2 \times S^1,x_0)$ via conjugation on both
summands, and the action of $\mathbb{Z}_2$ by permuting the two
summands. There is a natural isomorphism from $\pi_1(\mathbb{H}^2
\times S^1,x_0)$ to $\mathbb{Z}$ by mapping the fiber thought
$x_0$ to $1$. So we can identify $\mathcal{B}$ with the set of
unordered pairs of integers. Denote such a pair by $[a,b]$. Then,
$\mathcal{B}=\{[a,b]|~a,b\in\mathbb{Z}\}$. We denote by
$\mathbb{Z}[\mathcal{B}]$ the free $\mathbb{Z}$-module generated
by $\mathcal{B}$.

Now, let $K$ be a generic knot such that the only singularities of
its projection on $\mathbb{H}^2$ are transversal double points. We
equip $K$ with the canonical framing. Then the $I$ invariant of
the framed knot $K$ is defined by
\begin{eqnarray*}
I(K) & = & r(K)([h(K)+1,-1]-[h(K)-1,1]) \\
& & + \sum_d 2([h(K_d^+),h(K_d^-)-1] - [h(K_d^+)-1,h(K_d^-)]),
\end{eqnarray*}
where the sum is taken over all transversal double points of
$Pr(K)$. The invariant $I$ assigns an element of
$\mathbb{Z}[\mathcal{B}]$ to every framed knot in $\mathbb{H}^2
\times S^1$, and is invariant under isotopy of framed knots.

Now let $m:\mathbb{Z}[\mathcal{B}]\rightarrow\mathbb{Z}$ be the
$\mathbb{Z}$-linear mapping such that $m([a,b])=ab$ for $\forall
~[a,b]\in\mathcal{B}$. Then a simple calculation reveals the
following.

\begin{proposition}\label{w}
Let $K$ be a generic knot in $\mathbb{H}^2 \times S^1$ equipped
with the canonical framing. Then $m(I(K)) + 2w(K)+2r(K)h(K) = 0$.
\end{proposition}

With Proposition \ref{w} in hand, it's clear that Theorem
\ref{main} is a direct consequence of Proposition \ref{l} below.

\begin{proposition}\label{l}
Let $K$ be a generic knot in $\mathbb{H}^2 \times S^1$ equipped
with the canonical framing. Then $m(I(K)) + 2\beta(K) = 0$.
\end{proposition}

\begin{figure}[h]
  \includegraphics[width=4in]{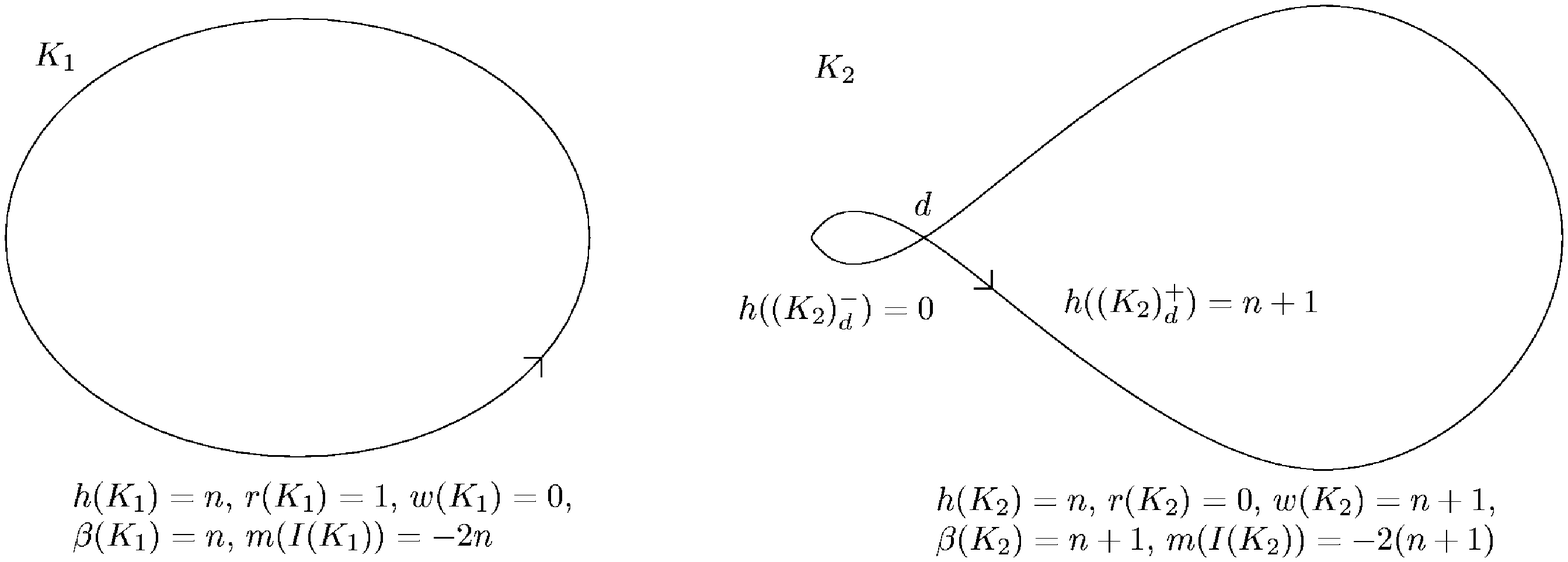}
  \caption{}\label{2knots}
\end{figure}

\begin{proof}
Note that $\beta(K)$ is indeed an invariant of framed knots, i.e.,
if two generic knots equipped with the canonical framing are
isotopic as framed knots, then the value of their
$\beta$-invariants are the equal. So, $m(I(K)) + 2\beta(K)$ is
also an invariant of framed knots. Let $K_s$ be a framed singular
knot in $\mathbb{H}^2 \times S^1$, whose only singularity is a
transversal self-intersection of two strands. Up to isotopy, there
are two ways to resolve this singularity. Let the two strands
intersecting at the singularity be $L_1$ and $L_2$. As in
\cite{Cher}, we denote by $K_s^+$ the resolution of $K_s$ such
that the tangent vector of $L_1$, the tangent vector of $L_2$, and
the vector from $L_2$ to $L_1$ form a position basis for
$\mathbb{H}^2 \times S^1$, and denote by $K_s^-$ the other
resolution of $K_s$. Then, $\beta(K_s^+)-\beta(K_s^-)=2$, and,
according to Theorem 1.2.5(2) of \cite{Cher},
$m(I(K_s^+))-m(I(K_s^-))=-4$. This shows that, if two framed knots
$K_1$ and $K_2$ are differed by a crossing of two transversal
strands, then $m(I(K_1)) + 2\beta(K_1)=m(I(K_2)) + 2\beta(K_2)$.
Since any homotopy between two framed knots, after a possible
slight perturbation, can be decomposed into a sequence of
isotopies and crossings of transversal strands, $m(I(K)) +
2\beta(K)$ is indeed a homotopic invariant of framed knots. Note
that any framed knot of homology $n$ is homotopic to one of the
two knots depicted in Figure \ref{2knots} equipped with canonical
framing. A straightforward calculation shows that $m(I(K)) +
2\beta(K) = 0$ for both of these knots. Thus, $m(I(K)) + 2\beta(K)
= 0$ for all framed knots in $\mathbb{H}^2 \times S^1$.
\end{proof}

\section{Computing $tb$ and $\mu$}

Since the standard hyperbolic metric on $\mathbb{H}^2$ has
curvature $-1<0$, its Levi-Citita connection form is a positive
contact form on $ST\mathbb{H}^2$. Under the identification $\Phi$
and $\Psi$, the corresponding contact structures on
$\mathbb{H}^2\times S^1$ and $\mathbb{R}_0^3$ are defined by
$\frac{dx}{y}+d\theta$, and $\frac{dz}{r^2}+d\varphi$,
respectively. Note that the latter is the restriction of $\xi_1$ on
$\mathbb{R}_0^3$.

Clearly, every Legendrian knot in $\mathbb{H}^2\times S^1$ is
generic. So, every invariant in this paper induces an invariant of
Legendrian knots in $\mathbb{H}^2\times S^1$. The Thurston-Bennequin number $tb$ and the Maslov
index $\mu$ in $(\mathbb{R}^3,\xi_1)$ restrict to invariants of Legendrian knots
in $\mathbb{H}^2\times S^1$. Their relations with the invariants
defined in Section \ref{inv-def} are described in Theorem \ref{contact}.

\begin{theorem}\label{contact}
Let $L$ be a Legendrian knot in $(\mathbb{R}_0^3,\xi_1|_{\mathbb{R}_0^3})$. Then,
\begin{enumerate}
    \item $tb(L)=\beta(L)=w(L)+r(L)h(L)$,
    \item $\mu(L)=r(L)+h(L)$.
\end{enumerate}
\end{theorem}

\begin{proof}
Note that $V$ is transverse to the contact plane field. So, by
definition, $tb(L)=\beta(L)$. The rest of part (1) is just Theorem
\ref{main}.

Let
$\pi:\mathbb{R}^3\rightarrow\mathbb{R}^2=\{(z,r,\varphi)|~z=0\}$
be the natural projection. Then $\mu(L)$ equals the degree of the
Gauss map of $\pi \circ L$. Define the map
$G:\mathbb{R}^3\rightarrow\mathbb{R}^3$ by
$G(z,r,\varphi)=(z,r^2,\varphi)$. Then $\pi \circ L$ and $\pi
\circ G \circ L$ are homotopic through immersed curves in
$\mathbb{R}^2$. So their Gauss maps have the same degree.

Let $(u,v,z)$ be the standard rectangular coordinates of
$\mathbb{R}^3$, i.e., $(u,v,z)=(r\cos{\varphi},r\sin{\varphi},z)$.
Then the map $\pi\circ G\circ\Psi:\mathbb{H}^2\times S^1
\rightarrow \mathbb{R}^2$ is given by
\[
\left\{%
\begin{array}{ll}
    u=y\cos{\theta},  \\
    v=y\sin{\theta}.  \\
\end{array}%
\right.
\]

Let $t$ be a regular parameter of $S^1$. We calculate the degree
$d$ of the Gauss map of $\pi \circ G \circ L$. First, we have
\[
d ~=~
\frac{1}{2\pi}\int_{S^1}\frac{u'v''-u''v'}{(u')^2+(v')^2}~dt.
\]
Then we substitute $(u,v)$ by $(y,\theta)$, and apply the relation
$x'+y\theta'=0$. This gives
\begin{eqnarray*}
d & = & \frac{1}{2\pi}\int_{S^1} \theta' ~dt +
\frac{1}{2\pi}\int_{S^1}
\frac{-yy''\theta'+(y')^2\theta'+yy'\theta''}{(y')^2+y^2(\theta')^2}
~dt
\\
& = & \frac{1}{2\pi}\int_{S^1} \theta' ~dt +
\frac{1}{2\pi}\int_{S^1} \frac{x'y''-y'x''}{(x')^2+(y')^2}~dt.
\end{eqnarray*}
But
\[
\mu(L) ~=~ d,
\]
\[
h(L) ~=~ \frac{1}{2\pi}\int_{S^1} \theta' ~dt,
\]
and
\[
r(L) ~=~ \frac{1}{2\pi}\int_{S^1}
\frac{x'y''-y'x''}{(x')^2+(y')^2}~dt.
\]
This proves part (2).
\end{proof}

\end{document}